\documentclass[12pt,a4paper]{article}

\usepackage{pdflscape}
\usepackage{amsmath}
\usepackage{amssymb}
\usepackage{latexsym}
\usepackage{srcltx}
\usepackage{graphics}
\usepackage{color}
\usepackage{epsfig}
\usepackage{color}
\usepackage{hyperref}

\newcommand{\re}{{\mathbb R}}
\newcommand{\n}{{\mathbb N}}

\newcommand{\z}{{\mathbb Z}}

\newcommand{\cT}{{\mathcal{T}}}

\newcommand{\cS}{{\mathcal{S}}}

\newcommand{\cU}{{\mathcal{U}}}

\newcommand{\bx}{{\boldsymbol{x}}}
\newcommand{\by}{{\boldsymbol{y}}}
\newcommand{\bz}{{\boldsymbol{z}}}
\newcommand{\be}{{\boldsymbol{e}}}

\newcommand{\ba}{{\boldsymbol{a}}}
\newcommand{\bb}{{\boldsymbol{b}}}
\newcommand{\bc}{{\boldsymbol{c}}}
\newcommand{\bv}{{\boldsymbol{v}}}
\newcommand{\bu}{{\boldsymbol{u}}}

\newcommand{\bn}{{\boldsymbol{n}}}
\newcommand{\bh}{{\boldsymbol{h}}}

\newcommand{\bss}{{\boldsymbol{s}}}

\newtheorem{theorem}{Theorem}
\newtheorem{prop}{Proposition}
\newtheorem{lemma}{Lemma}
\newtheorem{cor}{Corollary}
\newtheorem{remark}{Remark}

\newtheorem{defi}{Definition}

\date{}

\author{
Vladimir Yu. Protasov
\thanks{DISIM, University of L'Aquila, Italy;   Moscow State University, Department of Mechanics and Mathematics, Russia, {e-mail: \tt\small
vladimir.protasov@univaq.it}} , 
Tatyana Zaitseva
\thanks{Moscow Center for Fundamental and Applied Mathematics, Moscow State University, Department of Mechanics and Mathematics, Russia {e-mail: \tt\small
zaitsevatanja@gmail.com}} 
}

\title{Tiling of polyhedral sets
\thanks{
The first  author is supported  by the Russian
Foundation for Basic Research, projects  no. 19-04-01227 and 20-01-00469. 
}}

\begin{document}
\maketitle

\begin{abstract}

A self-affine tiling of a compact set~$G$ of positive Lebesgue measure is its partition to parallel shifts of a compact set which is affinely similar to~$G$.  We find all 
polyhedral sets (unions of finitely many convex polyhedra) that admit self-affine 
tilings. It is  shown that in~$\re^d$ there exist an infinite family of such 
polyhedral sets, not affinely equivalent to each other.   A special attention is paid 
to an important particular  case  when the matrix of affine similarity and the translation vectors are integer. Applications to the approximation theory and to the functional analysis are discussed.

\bigskip

\smallskip

\begin{flushright}

\end{flushright}

\end{abstract}
\bigskip

\vspace{1cm}

\begin{center}

\large{\textbf{1. Introduction}}	
\end{center}
\bigskip 

A {\em tiling} of a compact set~$G \subset \re^d$ is its partition to finely many 
sets that are disjoint (up to measure zero) and are translations of one compact set 
of positive Lebesgue measure called {\em tile}. If the tile is similar to~$G$ by means 
of some affine operator, then the tiling is {\em self-affine}.  

Self-affine  tilings have been studied in an extensive literature.  
Most of known tiles, apart from parallelepipeds, 
have fractal-like properties, which is natural due to their self-similarity. 
An important  problem is to find possibly simple sets that admit self-affine tilings.
For example, disc-like sets~\cite{BW}, polyhedral sets~\cite{GM, NM, YZ1, Zai}, etc.  
Besides  a geometric interest, this question has obvious applications 
in the space tiling, the crystallography,  the theory of functional equations, approximation theory, in constructing of orthonormal 
functional systems, in particular,  Haar systems and wavelets, etc. For example, if~$G$ possesses a piecewise-smooth boundary, then its characteristic function~$\chi_{G}$ 
has the maximal possible (among  piecewise-constant functions) regularity in~$L_2(\re^d)$: 
 its H\"older exponent is $\frac12$. 
Hence, the Haar system generated by~$\chi_{G}$  has the best approximation properties, the corresponding subdivision scheme and the cascade algorithm have the fastest rate of convergence, etc., see~\cite{KPS, Woj}. This holds, in particular, for polyhedral tiles. Not to mention, of course, that 
polyhedral sets are  convenient in practical applications. 

In this paper we classify polyhedral sets that admit self-affine tilings. 
The first results in this direction originated with Gr\"ochenig, Madych~\cite{GM},
who studied self-affine tilings of parallelepipeds and related Haar bases in~$L_2(\re^d)$. 
A complete classification of linear operators and parallel translations 
generating  self-affine tilings of parallelepipeds was done in~\cite[Theorems 1, 2]{Zai}. 
Other polyhedral tilings were addressed in~\cite{NM, YZ1, Zai}. 

This is relatively simple to show that among convex polyhedra, only parallelepiped 
admits a self-affine tiling. For non-convex polyhedra, 
the problem is more complicated. The two-dimensional case was done in~\cite[Theorem 3]{Zai} and the conclusion is the same: there are no self-affine polygonal tilings different from 
parallelograms. For higher dimensions, the corresponding conjecture was left open~\cite[Conjecture 1]{Zai}. 
Recently, in~\cite{YZ1} it was proved that if a set admits a self-affine tiling and has at least one convex polyhedral corner (conic intersection with a small ball) then this set is equivalent to a union of integer shifts of a unit cube. This strong result, however, does not solve the 
problem of characterising polyhedral self-affine tiles.   
Already in $\re^2$ there are polyhedral sets without convex corners. Moreover,  
in~$\re^3$ there are (non-convex) polyhedra without convex corners. For example, at each vertex of a regular tetrahedron 
we cut off a small tetrahedron and replace it by
 a dimple of the same form. We obtain a polyhedron with $16$ vertices none of which  has a convex corner. 
 
 On the other hand, there is a variety of disconnected polyhedral sets (unions of several convex polyhedra) that do admit self-affine tilings~\cite[Section 10]{Zai}. But their complete classification has been obtained only for one-dimensional case, see \cite{Long} which is based on \cite{Bruij}, \cite[Theorem 8]{Zai}.  
 In this paper we classify all polyhedral sets in~$\re^d$ with self-affine tilings (Theorem~\ref{th.10}). Moreover, we also classify the integer polyhedral attractors (Theorem~\ref{th.20}). This, in particular, gives  a complete characterization of Haar bases with polyhedral structure in~$L_2(\re^d)$. 
 
 It is worth mentioning that the the self-affine partitions of polyhedra  have also been studied under less restrictive conditions, when not only parallel translations but also 
 rotations are allowed~\cite{BV, BM, Sol, TD}. 
  
\bigskip

\begin{center}

\large{\textbf{2. The main results}}
\end{center}
\bigskip 

\begin{defi}\label{d.20}
Let~$G$ be a closed subset of~$\re^d$ of positive Lebesgue measure.  
A  {\em tiling} of~$G$ is its partition to a union 
of  compact sets $G \, = \, \cup_{i} T_i$,  such that  
all~$T_i$ are parallel shifts of each other and their  pairwise intersections are of Lebesgue measure zero. A tiling is {\em self-affine} if the sets~$T_i$ are affinely similar to~$G$.
\end{defi}
We consider only the  finite tilings, in which case~$G$ is compact. 
For example, all self-affine tilings are finite. 
A tiling will be denoted as $\cT = \{T_i\}_{i=1}^N$; each set~$T_i$
is called a {\em tile}.    
If the tiling is self-affine, then the similarity of each tile~$T$ to~$G$ 
is realized by an affine transform~$A: T\to G$. For all tiles,  
 those transforms have the same linear part $M$. 
 The matrix of the linear operator~$M$ (the dilation matrix) is denoted by the same symbol~$M$ and is supposed to 
 be expanding, i.e, all its eigenvalues are bigger than one in absolute value.   
 \bigskip

\begin{defi}\label{d.10}
A {\em convex polyhedron} is a compact subset of~$\re^d$ with a nonempty interior defined by several linear inequalities. A {\em polyhedral set} is a  union of finitely many convex polyhedra called {\em composing polyhedra}.  
\end{defi}

In one-dimensional case all polyhedral sets are unions of several segments.  
All such sets that admit self-affine tilings were characterised in~\cite{Zai}, it is also reduced to the result from \cite{Long}. 
To formulate this result we introduce some further notation. 
For natural numbers $a, n$, we denote 
$\cS (a, n) \, = \, \{k a \ | \ k = 0, \ldots , n-1\}$.
This is an    arithmetic progression of length~$n$ with the 
difference~$a$ starting at zero . For a natural~$r$ and positive vectors $\ba, \bn \in \z^r$, 
let $\cS (\ba, \bn)\, = \, \cS (a_1, n_1)+ \ldots + \cS (a_r, n_r)$
be the Minkowski sum of the progressions~$\cS (a_i, n_i)$. 
A pair of vectors~$\ba, \bn \in \z^r$  are called {\em admissible}
if  $a_1 = n_1= 1$ and for each $i \ge 2$, we have 
$a_i\ge 2, n_i \ge 2$, and   $a_i$ is divisible by~$a_{i-1}n_{i-1}$. 
Now we formulate the characterization of univariate polyhedral tilings~\cite[Theorem~8]{Zai}:  
 \smallskip 
 
{\em A union of finitely many segments in~$\re$ denoted by~$G$  possesses  a self-affine tiling 
if and only if there are $r\in \n$ and an admissible pair of 
vectors~$\ba, \bn \in \z^r$ such that  $G$  is equivalent (up to normalisation) to the set: 
\begin{equation}\label{eq.one-dim}
\Bigl\{  \  [k, k+1], \quad  k \, \in \, \cS(\ba, \bn)\, \Bigr\}\, . 
\end{equation} }
\smallskip

Thus, a polyhedral subset of~$\re$ that admits a self-affine tiling is 
equivalent to a disjoint union of integer shifts of the unit segment.
The number of segments is~$n_1\ldots n_r$. 
The case of one segment corresponds to $r=1$.  
\smallskip

Our first  result classifies all polyhedral sets in~$\re^d$ that possess 
self-affine tilings. 
\begin{theorem}\label{th.10}
A polyhedral set in~$\re^d$  possesses  a self-affine tiling 
if and only if it is affinely equivalent to a disjoint union of integer translates of the unit cube which is a direct product of~$d$ 
one-dimensional sets of the form~(\ref{eq.one-dim}).  Each of those $d$ sets corresponds to its own triple~$(r, \ba, \bn)$, where $\ba, \bn$ is an admissible pair from~$\z^r$.  
\end{theorem}

\begin{figure}[ht]
\centering
\includegraphics[width = 0.45\textwidth]{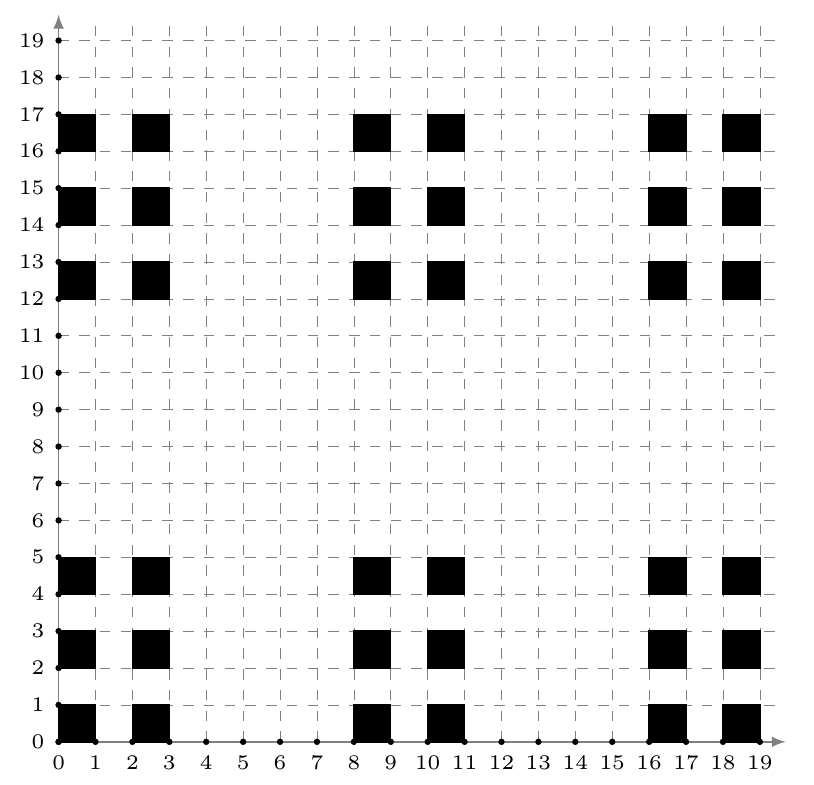} 
\caption{The two-dimensional example of a polyhedral set.}
\label{polyhedra2}
\end{figure}

The example of a polyhedral set in $\re^2$ is given on Fig. \ref{polyhedra2}. The first triple (axis $x$) is (3; (1, 2, 8); (1, 2, 3)), the second one (axis $y$) is (3; (1, 2, 12); (1, 3, 2))). 

The second result characterises polyhedral sets with special self-affine tilings. 
\begin{defi}\label{d.30}
A self-affine tiling~$\cT \, = \, \{T_i\}_{i=1}^m$ is called {\em integer} if the dilation matrix $M$ has integer entries and $G = MT_i + \bss_i$, where $m = |{\rm det}\, M|$ and 
  $\bss_1, \ldots , \bss_m$ are integer vectors representing  all different 
  quotient classes in~$\z^d / M\z^d$. Those vectors are called {\em digits} and 
  compose the set of digits  $D \, = \, \{\bss_i\}_{i=1}^m$. The set~$G$ that possesses an integer tiling is called an {\em integer attractor}. 
\end{defi}
Thus, in the integer tiling the matrix~$M$ is integer and 
the shift vectors compose a complete set of digits with respect to~$M$. 
This means that $\bss_i - \bss_j \notin M\z^d$ for all~$i\ne j$ and for every~$\bss \in \z^d$
there exists~$j$ such that $\bss- \bss_j \in M\z^d$. Integer attractors are applied in the combinatorics and number theory as well as in the constuction of orthenormal Haare beses 
in~$L_2(\re^d)$. 
By Theorem~\ref{th.10}, a polyhedral integer attractor must be equivalent to 
a union of integer shifts of the unit cube. However, among all the sets described 
in Theorem~\ref{th.10} only one case corresponds to an integer tiling. 
The following theorem gives a complete classification.  
\smallskip

\begin{theorem}\label{th.20}
If a polyhedral set is an integer attractor, when it is a parallelepiped. 
\end{theorem}
Note that one integer attractor can be generated by different dilation matrices and by sets of digits. 
\bigskip 

\begin{center}

\large{\textbf{3. The roadmap of the proofs}	}
\end{center}
\bigskip 

The proof of Theorem~\ref{th.10} consists of several steps. 
We shall briefly describe them below and then in Section~5 give a detailed proof. The proof of Theorem~\ref{th.20} is much shorter (of course, after referring to Theorem~\ref{th.10}); its idea is described in the end of this section. A complete proof of Theorem~\ref{th.20} is given in Section~6.  
\bigskip 

\textbf{The steps of  the proof of Theorem~\ref{th.10}}. Let us have a polyhedral set~$G \subset \re^d$ with a self-affine tiling~$\cT$. 
\bigskip 

\noindent {\tt Step~1}. The convex hull of the set~$G$  is a convex 
polyhedron. Vertices  of this polyhedron are referred to as  {\em extreme vertices}
of~$G$. For an    arbitrary extreme vertex~$\bv$, we 
consider the corresponding {\em corner}~$K$ of~$G$, which is the intersection 
of~$G$ with a small ball centered at~$\bv$. This is a cone with an apex~$\bv$, 
maybe non-convex. The set ${\rm co}\, (K)$ is a corner of~$\, {\rm co}\, (G)$; 
it is convex and has an apex at the vertex~$\bv$. In the partial order defined by the 
cone~${\rm co}\, (K)$ (we keep the same notation for a corner and for the corresponding cone) the vertex~$\bv$ is a unique minimal element of the polyhedron~${\rm co}\, (G)$. 
This fact implies that there exists a unique tile
$T \in \cT$ that contains~$\bv$ (see Proposition~\ref{p.30}). 
Moreover, a sufficiently small corner of~$G$ at the extreme vertex~$\bv$  is also a corner of~$T$, and 
$\bv$ is an extreme vertex for~$T$ as well
 (Proposition~\ref{p.40}).  
 \medskip 

\noindent {\tt Step~2}. Thus, for every extreme 
vertex~$\bv$ of~$G$, there exists a unique tile from~$\cT$  
containing~$\bv$, for which $\bv$ is also an extreme vertex with the same corner. This defines a map from the set of extreme vertices of $G$ to extreme vertices of~$T$. Since the tiling is self-affine, $T$ is affinely similar to~$G$ and hence this 
map defines a map of the set of extreme vertices of~$G$ to itself. 
 The graph corresponding to this map has one outgoing edge from each vertex and hence has a cycle of some length~$n$. Therefore,
for the $n$th  iteration of the tiling~$\cT$, there 
  exists a {\em stationary} vertex~$\bv$  corresponding to itself. 
Hence, it corresponds to the same vertex of~$T$. Moreover, after extra iterations of the tiling, it
  may be assumed that the convex corner ${\rm co}\, (K)$
   has the same corresponding faces in~$G$ as in~$T$. 
   This is proved in 
 Proposition~\ref{p.45}. 

 \medskip 

\noindent {\tt Step~3}. Thus, $G$ has at least one {\em stationary vertex}, which corresponds to the same extreme vertex of $G$ and of~$T$. Our goal is to prove that the corner~$K$ at this vertex is convex and simple (has exactly~$d$ edges). To this end, we first show that if $K$ contains a face of the convex corner~${\rm co}\, (K)$, then the intersection of all  tiles with that face define a self-affine tiling on it 
(Lemma~\ref{l.10}). This will allow us to apply induction arguments in the dimension of the faces. 

 \medskip 

\noindent {\tt Step~4}. If the corner $K$
at a stationary vertex~$\bv$ contains a $j$-dimensional face 
of~${\rm co}\, (K)$, then it contains a parallelepiped 
in~$L$ with the same corner at~$\bv$  (Lemma~\ref{l.20}). 
 \medskip 

\begin{figure}[ht]
\centering
\includegraphics[width = 0.36\textwidth]{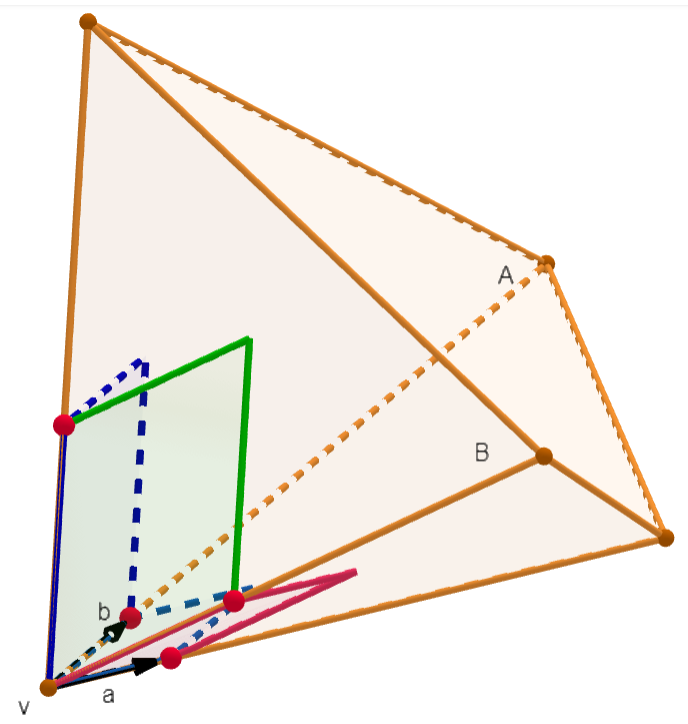} 
\caption{The ``lily'': the parallelepipeds on facets spanned by vectors along the edges.}
\label{lili1}
\end{figure}

\noindent {\tt Step~5}. This is a  keystone of the proof. We 
establish the following auxiliary geometrical fact, which is, probably, of some independent interest. 
For a convex polyhedral cone~$C$ with simple facets and for 
arbitrary vectors going from the apex along the edges  (one vector on each edge), the following holds: either $C$ is simple, 
or there exist two facets $A, B$ of $C$
and two vectors~$\ba, \bb$ from our family on different edges such that  the shifted parallelepipeds 
$\ba + P(A)$ and $\bb + P(B)$ have a common interior point, where $P(X)$ is the parallelepiped spanned by the vectors on a facet~$X$ (see Fig. \ref{lili1}).  In the former case ($C$ is simple) those parallelepipeds form a boundary of a full-dimensional parallelepiped with the corner~$C$. This is Lemma~\ref{l.30}, which  can be called ``the lily lemma''
since the parallelepipeds form  ``lily petals''. 

 \medskip 

\noindent {\tt Step~6}. Applying the result of the previous step we prove by induction
in~$j$ the following statement: if a vertex~$\bv$
is stationary, then for each $j  = 1, \ldots , d$, the corner~$K$ contains all faces of dimension~$j$ of the convex corner~${\rm co} \,(K)$ and all those faces are simple (Proposition~\ref{p.50}). 
Taking $j=d$ we immediately obtain that the corner~$K$ is convex and simple. 
 \medskip 
 
\noindent {\tt Step~7}. Thus, the corner of $G$ at every stationary vertex is convex and simple. Hence, by \cite[Theorem~1.9]{YZ1}, 
 the set $G$ is equivalent to a union of several integer shifts of the unit cube. Then we establish Lemma~\ref{l.40} according to which there exists a subset of the tiling~$\cT$
 that forms a tiling  of a rectangular parallelepiped. 
 Then we apply the main result of the paper \cite{Nath} on discrete tiling of a cube and conclude that 
 $T$ is a direct product of $d$ one-dimensional tiles. Note that in case of infinite discrete tilings of the set $\n \times \n$ this does not hold, see \cite{Niv}. 
 Hence so is~$G$. It remains to invoke the classification of 
 one-dimensional tilings from~\cite[Theorem 8]{Zai} formulated in Section 2.

\bigskip 

\textbf{The idea of  the proof of Theorem~\ref{th.20}}.
Applying Theorem~\ref{th.10} we obtain that $G$
 is affinely similar to a direct product of special sets of the form~(\ref{eq.one-dim}). 
 In particular, this yields that 
 $G$ is equivalent to a union of disjoint integer shifts of the unit cube. 
  If the tiling~$\cT$
 is integer, then all its iterations are also integer. Taking a sufficiently big iteration 
 we can assume that all tiles have diameter less than one. Hence, each tile is contained in a unique cube. This implies that 
 the shifts cannot be from different quotient classes, provided $G$ contains at least two cubes.
 The proof is in Section~6.  
\bigskip 

\begin{center}

\large{\textbf{4. Notation and preliminary facts}	}
\end{center}
\bigskip

We use convex separation theorems. 
Let $X$ and $Y$ be subsets of~$\re^d$. The set 
 $X$ is separated from~$Y$ by a nonzero element  $\bc \in \re^d$ if 
$(\bc, \bx) \le (\bc, \by)$ for all $\bx \in X, \, \by \in Y$. 
This separation is {\em strong} if the set $Y\setminus X$ is nonempty and  
$(\bc, \bx) < (\bc, \by)$ for all $\bx \in X, \, \by \in Y\setminus X$. 
By the convex separation theorem, if $X, Y$ are both convex, 
the interior of $Y$  is non-empty and  
does not intersect~$X$, then $X$
 can be separated from~$Y$. Every face of a convex polyhedron~$G$ is strongly separated from it. 
A point on the surface of~$G$ is strongly separated from it
precisely when it is extreme. 

\newpage

 \bigskip 

\begin{center}
\textbf{4.1. Cones}	
\end{center}
\bigskip

A {\em cone}~$K$ with the apex at the origin is a closed subset of $\re^d$ such that if $\bx \in K$
then for each $\lambda \ge 0$, we have $\lambda \bx \in K$. A cone is nondegenerate if 
it possesses a nonempty interior. A cone is {\em pointed} if it does not contain a 
straight line. A cone is {\em convex} if for every~$\bx, \by \in K$,  we  have 
$\bx + \by \in K$. Any ray that belongs to the boundary of a cone is called its {\em generatrix}. 
If a generatrix does not belong to a convex hull of  other generatrices, it is called 
{\em an extreme edge}. In particular, a generatrix of a convex cone is an extreme edge if it does not belong to a linear span of two other generatrices. In this case we drop the word ``extreme'' and  call it edge.

 In what follows we always assume that a convex cone is nondegenerate and pointed. In this case it can be strongly separated from its apex.

A convex cone~$K\subset \re^d$ is called polyhedral if it is defined by a system of 
linear inequalities. It has faces of all dimensions from zero 
(the apex) to $d$ (the whole cone).  A {\em facet} is a face of dimension~$d-1$; an edge is a face of dimension one. A cone  is {\em simple} if is is affinely equivalent to 
the positive orthant~$\re^d_+$. So, a simple cone has exactly $d$ facets and $d$ edges. A cone possesses simple facets if all their facets  are simple 
$(d-1)$-dimensional cones. Clearly, in this case its  faces of all dimensions $\le d-1$ are simple. 

A  {\em corner}  is an intersection of a cone  with a  ball centered at the apex. The cone is an {\em extension} of 
every its corner.     
 We often identify a corner and its extension and use the same notation for them.

Every convex cone defines a partial order in~$\re^d$ as follows: 
$\bx \ge \by$ if $\bx - \by \in K$. In particular, $\bx \ge 0$ if $\bx \in K$. 
Since $K$ is pointed, it follows that if $\bx \ge \by$, then $\by \ge \bx$ is impossible, 
unless $\bx = \by$. A point $\bx$ is called {\em minimal} for a set $\Omega \in \re^d$
if $\bx \le \by$ for all~$\by \in \Omega$. Not every set possesses a minimal element, but 
if it exists, it is unique. Indeed, if there are two minimal elements $\bx$ and $\by$, 
then $\bx \le \by$ and $\by \le \bx$, hence $\bx = \by$.   

Let $B$ be a subset of a hyperplane~$L$ with a nonempty (in~$L$) interior and 
$\bv\notin L$ be a point. Then all segments 
connecting $\bv$ with points from~$B$ form a {\em bounded cone}
with the apex~$\bv$ and a base~$B$.  

\begin{center}

\textbf{4.2. Polyhedral sets}	
\end{center}
\bigskip

A  {\em corner} of a polyhedral set is its intersection  with a small ball centered at a vertex~$\bv$ of some composing polyhedron. 
We always assume that the ball is small enough and 
intersects only those faces of the composing polyhedra adjacent to~$\bv$.  
Clearly, the extension of a corner is nondegenerate, but possibly  
non-convex and not pointed. 

\begin{figure}[ht]
\centering
\includegraphics[width = 0.45\textwidth]{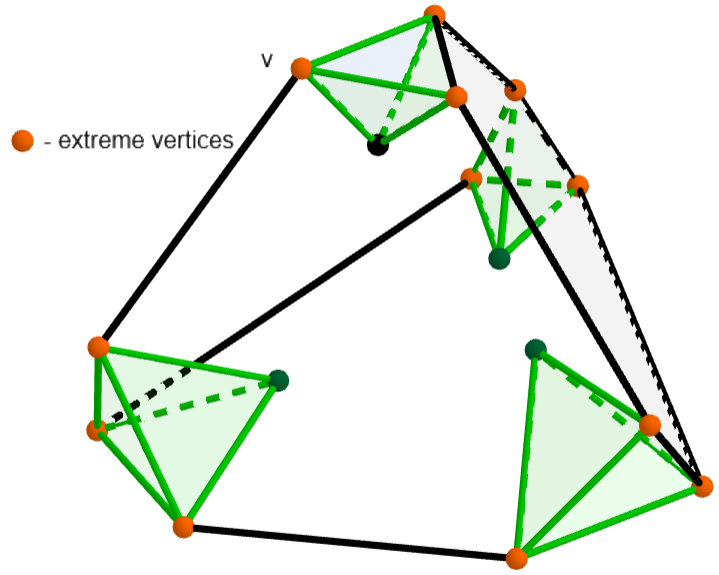} 
\caption{An example of extreme vertices of modified tetrahedron without convex corners.} 
\label{extreme}
\end{figure}

A point of a polyhedral set~$G$  
is said to be  {\em an extreme vertex} if it is a vertex 
of~${\rm co }\, (G)$ (see Fig. \ref{extreme}). 
A point $\bv \in G$ is an extreme vertex if and only if it  can be 
strongly separated from~$G$.   
That vertex is called {\em convex} if it is a vertex of a convex corner of~$G$
(as usual, a convex corner is assumed to be non-degenerate and pointed).  

Similarly we define a {\em composite extreme face}~$L$ of a polyhedral set~$G\subset \re^d$. 
This is  the intersection of $G$ with a separating hyperplane of~$L$. 
For the sake of simplicity we usually call~$L$ {\em extreme face}. 
It is a union of several convex sets that are faces of polyhedra 
that form~$G$. The maximal dimension of those sets is the 
{\em dimension of the extreme face}.  
An extreme face of dimension~$d-1$ is an {\em 
extreme facet}. Let~$G' \subset G$ be an extreme facet and $S$ be a convex subset of~$G'$; then 
  a {\em layer} of~$S$  is a bounded 
cone with   the base~$S$ that is contained in~$G$. Not every convex subset of a facet has a layer. 
For example, if $ABC$ is a triangle and 
$A', B', C'$ are the midpoints of its sides, then the union 
of the triangles $AB'C'$ and $BC'A'$ is a polyhedral set 
with the facet~$AB$. This facet, however, does not have a layer.

As we have already mentioned, there are polyhedra without convex vertices. However, 
extreme vertices always  exist. 
\begin{prop}\label{p.10}
If $K$ is a corner of a polyhedral set at its extreme vertex, then ${\rm co}\, (K)$
is a nondegenerate convex pointed cone.  
\end{prop}
{\tt Proof.} The convexity and nondegeneracy are 
obvious. The pointedness of ${\rm co}\, (K)$ follows from the 
pointedness of the corner of 
${\rm co}\, (G)$, which contains~${\rm co}\, (K)$. 

{\hfill $\Box$}
\begin{prop}\label{p.20}
Every polyhedral set  is contained in the convex hull of its extreme vertices.  
\end{prop}
{\tt Proof.}   Extreme vertices of a polyhedral set~$G$ are 
 vertices  of the convex polytope~${\rm co} (G)$. By the Minkowski theorem, 
 a convex polyhedron 
is a convex hull of its vertices. 

{\hfill $\Box$}
\begin{cor}\label{c.20}
Every polyhedral set in~$\re^d$ possesses at least $d+1$ extreme vertices.   
\end{cor}
\bigskip

 \begin{center}
\large{\textbf{5. Proof of Theorem~\ref{th.10}}	}
\end{center}
\bigskip 
 
  We begin with several basic facts, which 
  hold for all tilings, not necessarily self-affine.  
\begin{prop}\label{p.30}
Suppose $\cT$ is a tiling of a polyhedral set~$G$;
 then  every extreme vertex of~$G$ is contained in a unique element of~$\cT$. 
\end{prop}
{\tt Proof.}  Let an element $T \in \cT$ contain an extreme  vertex $\bv$ of~$G$. 
Denote by $\tilde K$ the extension of the corner of ${\rm co} (G)$ at the 
vertex~$\bv$.
Since ${\rm co} (G) \, \subset \, \tilde K$ and 
$T \subset G$, it follows that $T \subset \tilde K$. 
Hence,   $\bv$  is the minimal element of~$T$ in the order defined by the cone  
$\tilde K$. 
If another element $T'\in \cT$ also contains~$\bv$, then $\bv$ is the minimal element of~$T'$. 
From the uniqueness of the minimal element and from the fact that $T'$ is a parallel shift of 
$T$ by some nonzero vector~$\ba$, we see that $\bv+ \ba = \bv$;
 therefore, $\ba = 0$ and $T' = T$. 

{\hfill $\Box$}
\medskip 

\begin{prop}\label{p.40}
If an element $T \in \cT$ 
contains an extreme vertex~$\bv$ of~$G$, then it contains a sufficiently small corner at~$\bv$.
In particular, if the tiling is polyhedral, then $T$ and $G$ have the same corner at~$\bv$
and this vertex is extreme for~$T$.  
\end{prop}
{\tt Proof.}   By Proposition~\ref{p.30}, $T$ is a unique tile containing~$\bv$, hence 
all other elements of~$\cT$ are located on positive distances from~$\bv$. 
Choosing the radius of the ball smaller than all those distances, we obtain a 
corner that does not intersect other elements of~$\cT$. 
On the other hand, the tiling covers the whole set~$G$, therefore 
the intersection of the small ball with $G$ coincides with its  intersection 
with~$T$. Thus, $T$ and $G$ have the same corner at~$\bv$. 
Furthermore, since $\bv$ can be strongly separated from~$G$ by a hyperplane, 
the same hyperplane separates $\bv$ from~$T$, because~$T \subset G$. 
Consequently,~$\bv$ is an extreme vertex of~$T$. 

{\hfill $\Box$}
\medskip

Thus, every extreme vertex of~$G$ is associated to an 
extreme vertex of~$T$, which is the minimal element of~$T$ in the order defined by 
the corner of~${\rm co} (G)$ at that vertex. 
The inverse correspondence may not be well-defined: an extreme vertex of~$T$ can have no corresponding vertices~from~$G$ or can have several ones.

Now we turn to self-affine tilings. Each tile~$T$ is similar to~$G$
by means of some affine transform~$A: T \to G$. 
It maps each extreme vertex of~$T$ to a corresponding vertex of~$G$. For all tiles, those transforms have 
the same linear part defined by the dilation matrix~$M$. 
\begin{defi}\label{d.40}
Let a polyhedral set~$G$ possess a self-affine tiling~$\cT$.  
Then its extreme vertex~$\bv$ is called {\em stationary}  if   a unique tile~$T\in \cT$ containing $\bv$ possesses the following properties: 
\smallskip 

1) $\bv$ is a fixed point of the affine transform~$A: T\to G$; 

2) the transform~$A$ respects the corner~$K$ of~$G$ at the vertex~$\bv$  
and  respects  all faces of the cone ${\rm co} \,(K)$; 

3) the extension of~$K$ contains $G$. 
\end{defi}
Extreme vertices depend only on the polyhedral set~$G$, 
while the stationary vertices depend also on its self-affine tiling. 
As we know, $G$ has at least $d+1$ extreme vertices. However, for some tilings, none of those vertices are stationary. 
Nevertheless, at least one of them does become stationary
after certain iteration of the tiling. This is guaranteed by the following 
\begin{prop}\label{p.45}
For every self-affine tiling $\cT$ of a polyhedral set~$G$, 
there is $n\in \n$ such that~$G$ with the  tiling~$\cT^n$ 
possesses at least one stationary vertex.
\end{prop}
{\tt Proof.} Proposition~\ref{p.40} implies that every extreme vertex
$\bu$ of~$G$ corresponds to a unique extreme vertex~$\bu'$ of~$T$.
The latter, in turn, corresponds to a unique extreme vertex~$A\bu'$ of~$G$. 
Thus, we have a map~$\bu \mapsto A\bu'$ defined on the set of extreme vertices of~$G$. 
Denote this map by~$\varphi$. In general,~$\varphi$ may not be injective. 
 Iterating~$\varphi$, we obtain an extreme  vertex~$\bv$ of~$G$
and a number~$k$ such that $\varphi^k(\bv) = \bv$. This means that, for the 
tiling~$\cT^k$, the vertex~$\bv$ is covered by the corresponding 
(in the sense of similarity) vertex of the element of the tiling. 
Thus, after taking some power of the tiling 
we assume (keeping  the previous notation for the new tiling) 
 that $G$ and $T$ have the same corresponding vertex~$\bv$, i.e., $A\bv = \bv$. 
By Proposition~\ref{p.40}, the corners of $G$ and of $T$
at~$\bv$ coincide. 
Hence, the transform~$A$  preserves this corner~$K$ and  
therefore,  defines a permutation of its extreme edges. 
Some power of this permutation is identical. This means that passing 
again to a power of the tiling we may assume that 
$A$ maps each extreme edge of the cone~$K$ to itself. 
Hence, it respects all faces of~${\rm co}\, (K)$ and so, $\bv$ is a stationary vertex. 

With possible  further iteration of the tiling it may be assumed that the tile~$T$ is small 
enough and 
is contained in a small ball defining the corner~$K$. Hence, $T$ is contained in 
the extension of~$K$. Since $A$ maps $T$ to $G$ and respects $K$ it follows that, 
$G$ is also contained in the extension of~$K$. 

{\hfill $\Box$}
\medskip

In what follows we simplify the notation and assume that we are already given the 
$n$th power of the tiling. Thus, for a tiling~$\cT$ of the set~$G$, there exists a stationary vertex~$\bv$.


Until now we dealt with tiles covering the extreme vertices. 
Now we go further and  look at  extreme faces.
 Let us stress that 
 we consider faces of convex cones only. 
Let $K$ be a corner of a polyhedral set~$G$ 
and $L$ be a $j$-dimensional face of the convex cone~${\rm co} (K),\ j\, \ge \, 1$. 
We denote $G_L  =  G \cap L$.  If  $G$ admits a self-affine tiling, then 
$G$ lies in the extension of~$K$. Hence, in this case ${\rm co}\, (G_L)$ is 
a $j$-dimensional face of~${\rm co}\,(G)$.  
The following lemma reduces the dimension in the  proof of Theorem~\ref{th.10}. 
\begin{lemma}\label{l.10}
Let $\cT$ be a self-affine tiling of~$G$ and $T$ be its element containing
a stationary vertex~$\bv$. 
Suppose a corner~$K$ at~$\bv$ contains a face~$L$ of~${\rm co}\,(K)$; 
then all elements of $\cT$ intersecting~$G_L$ are translations of~$T$
by vectors parallel to~$L$. The sets of intersection form a self-affine tiling of~$G_L$. If  $T\cap L$ is a composite facet of $T$, then every convex subset of it has a layer in~$T$.   
\end{lemma}
{\tt Proof} in the Appendix.  In particular, for $j=1$, we obtain: 
\begin{cor}\label{c.30}
Let $\ell$ be an extreme edge of~$G$ going from a stationary vertex~$\bv\in G$.
Let $\cT$ be a self-affine tiling of~$G$ and $T$ be its element containing~$\bv$. 
Then all elements of $\cT$ intersecting~$\ell$ are translations of~$T$
by vectors parallel to~$\ell$. The sets of intersection form a 
tiling of~$\ell$.  
\end{cor}
{\tt Proof.} As we know, all edges of~${\rm co}(K)$
(faces of dimension one) are extreme edges of~$K$
and they lie in~$K$. Hence, the assumptions of  Lemma~\ref{l.10} for $\ell$ are fulfilled. 
Now applying Lemma~\ref{l.10} in the case~$j=1$ we conclude the proof. 

{\hfill $\Box$}
\medskip 

Now we  show that if a corner of~$K$ of a stationary vertex 
contains a face of the corner ${\rm co}\, (K)$, then the set~$T\cap L$
contains a special parallelepiped. 
\begin{lemma}\label{l.20}
Under the assumptions of Lemma~\ref{l.10}, suppose that 
$L$ is a simple $j$-dimensional cone with edges~$\ell_s, \, s = 1, \ldots , j$. 
Let~$\bb_s$ be the most distant from~$\bv$ point of the edge~$\ell_s$
for which 
the segment~$[\bv, \bb_s]$ is  contained in
$T$. Then $T$ contains a $j$-dimensional 
parallelepiped spanned by the segments $[\bv, \bb_s] \, , \ s = 1, \ldots , j$. 
\end{lemma}
{\tt Proof} is by the induction in the dimension~$j$. 
For $j=1$, the statement is obvious.  Assume it holds for some dimension~$n < j$. 
Consider the face $L_{n+1}$ spanned by the edges~$\ell_1, \ldots , \ell_{n+1}$
and its face~$L_n$ spanned by the first $n$ edges. 
Since $L_{n+1}$ is a simple cone, we can identify it with~$\re^{n+1}_+$
and assume that all the segments~$[\bv, \bb_s]$ are of length one. 
Denote $\bb_s - \bv = \be_s$. 
We use the same notation for the tiles from~$\cT$ and for their intersections with 
the face~$L_{n+1}$. 
We need to show that the unit cube~$P_{n+1} = \{\bx \in \re^{n+1}_+ \, , \, x_i \le 1, 
i = 1, \ldots , n+1\}$ is contained in~$T$. 
By the inductive assumption, $T$ contains its face~$P_{n} = 
\{\bx \in P_{n+1} \, , \, x_{n+1} =0\}$. Since by Corollary~\ref{c.30}, 
the shifts of the set~$T\cap \ell_{n+1}$ form a tiling of~$\ell_{n+1}$
it follows that the set~$T' = T+ \be_{n+1}$
is also a tile from~$\cT$. 

\begin{figure}[ht]
\centering
\includegraphics[width = 0.37\textwidth]{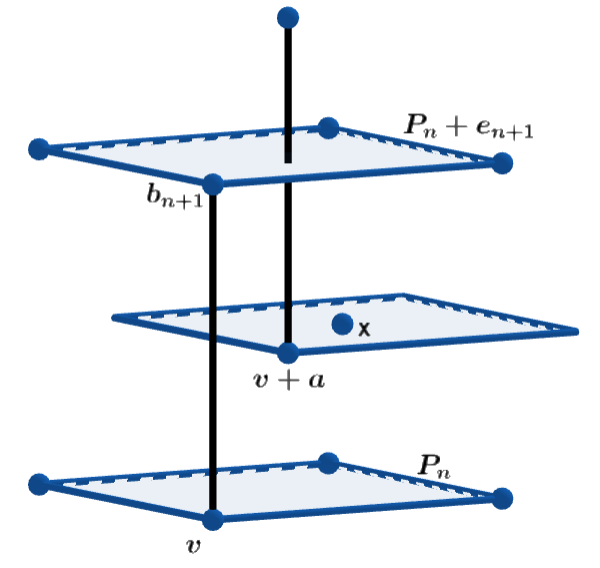} 
\caption{Illustration to the proof of Lemma 2.} 
\label{lemma2}
\end{figure}

If $P_{n+1}$ does not lie in~$T$, then there is a point $\bx \in  {\rm int}\, P_{n+1}$ which 
is not in~$T$ (see Fig. \ref{lemma2}). 
On the other hand, it must belong to 
some tile~$T'' =  T + \ba, \, \ba \in {\rm int}\, P_{n+1}$.  
For every $\by \in T'$, we have $y_{n+1} \ge 1$, hence   
none of the points $\bx, \bv + \ba$ belongs to~$T'$,  and so $T' \ne T''$.   On the other hand, $T''$ must contain the parallel shift of the segment  $[\bv, \bb_{n+1}]$ 
by the vector~$\ba$. This segment $[\bv + \ba, \bb_{n+1} + \ba]$  is of length one and it intersects 
the face $P_{n} + \be_{n+1}$ of $P_{n+1}$. 
However, this face  lies in~$T'$ with some layer (Lemma~\ref{l.10}). Hence, $T$ and $T'$ have a common interior point.  
The contradiction proves that $P_{n+1}$ is in~$T$, which completes  the inductive step.

{\hfill $\Box$}
\medskip 

\begin{figure}[ht]
\begin{minipage}[h!]{0.36\linewidth}
\centering
\includegraphics[width = 1\textwidth]{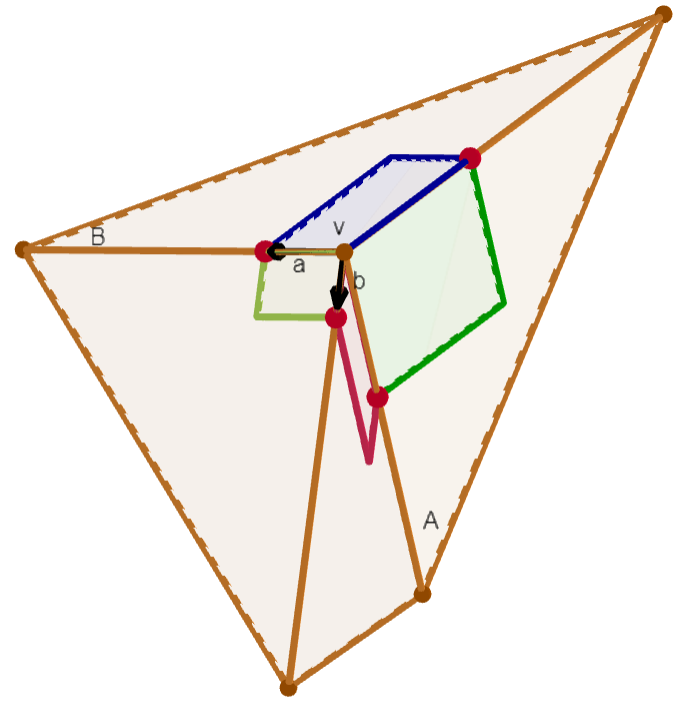} 
\caption{Parallelepipeds on facets spanned by the vectors on the edges.} 
\label{lili2}
\end{minipage}
\begin{minipage}[ht]{0.1\linewidth}
\centering
\quad
\end{minipage}
\begin{minipage}[h!]{0.36\linewidth}
\centering
\includegraphics[width = 1\textwidth]{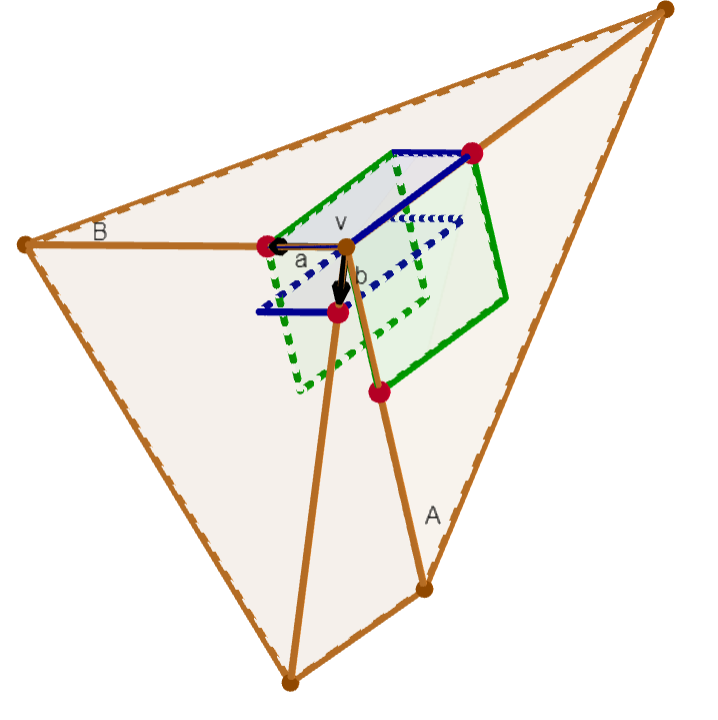}
\caption{Intersection of shifted parallelepipeds.} 
\label{lili3}
\end{minipage}
\end{figure}

 Now we are going to prove Proposition~\ref{p.50}, from which 
 Theorem~\ref{th.10} simply follows. We need the following geometrical lemma, 
 which is, probably,  of some independent interest. 
Let $C$ be a convex polyhedral cone with simple facets and with the apex at the origin~$O$. 
On every edge of~$C$  one chooses  an arbitrary  point $\bc \ne 0$,  
which is referred to as a {\em directing point} and the segment $[O, \bc]$
is a {\em directing} segment of that edge.  
For a given facet~$H$ of $C$, we denote by~$P(H)$ its 
{\em directing parallelepiped}, which is spanned by the 
directing segments of that facet. Thus, $P(H)$ is a $(d-1)$-dimensional 
parallelepiped contained in~$H$. The family of all such parallelepipeds
form a ''lily'' based on the cone~$C$.   
\begin{lemma}\label{l.30} (Lily lemma). 
For a convex polyhedral cone~$C$ with simple facets and for 
arbitrary directing points on its edges, the following holds: either $C$ is simple, 
or there exist two facets $A, B$
and two directing points $\ba, \bb$ on different edges of~$C$ such that  the shifted parallelepipeds 
$\ba + P(A)$ and $\bb + P(B)$ have a common interior point (see Fig. \ref{lili2}, \ref{lili3}). 
\end{lemma}

{\tt Proof} is in the Appendix.   Now we formulate the main proposition. 
\begin{prop}\label{p.50}
If a polyhedral set admits a self-affine tiling, then its corners at 
all stationary vertices are  convex and simple.
\end{prop}
 {\tt Proof.} 
 As usual,  $G$ denotes a polyhedral set, ~$\cT$ is a self-affine tiling of~$G$, 
 $K$ is a corner of~$G$
at a stationary vertex~$\bv$, and 
$T\in \cT$ is  the tile
covering~$\bv$. 
After possible iteration of the tiling it can be assumed that 
  the size of 
 $T$ is as small as needed. We prove the following statement which immediately implies Proposition~\ref{p.50}: 
 \smallskip 
  
 {\em For each $j  = 1, \ldots , d$, the corner~$K$ contains all faces of dimension~$j$ of the convex corner~${\rm co} \,(K)$ and all those faces are simple.}
  \smallskip 
  
  Applying this statement for $j=d$ we obtain that $K$ contains the whole convex cone~${\rm co} \,(K)$, 
  which is, moreover, simple. 
 So $K = {\rm co} \,(K)$, and hence $K$ is convex and simple.  
  The proof of this statement is by induction in the dimension~$j$. 
  \smallskip 
  
  $\mathbf{j=1}$. In this case the faces are the extreme edges of~${\rm co}\, (K)$, which, as we know, 
  coincide with the extreme  edges of~$K$. So, in this case the statement is true. 
   
   \smallskip 
  
  $\mathbf{j\, \to \, j+1}$. Assume the statement holds for some $j\le d-1$. Suppose the converse: the corner ${\rm co} \,(K)$ possesses  a face~$L$ of dimension~$j+1$ 
  that does not lie in~$K$.  Then the set~$S = L\setminus K$ is nonempty. 
  Since ${\rm co} \,(K)$ is small enough, it follows that $S$ is a corner 
  with apex~$\bv$. 
  
  Let $H$ be an arbitrary $j$-dimensional face of~$L$. It is also a face 
  of~${\rm co}\,(G)$ and by the inductive assumption, $H$ is simple and is contained in~$K$. 
For an arbitrary edge  of~$L$, we take its directing point~$\ba$, 
for which the segment~$[\bv, \ba]$ (the directing segment) is the biggest by inclusion segment 
of that edge contained in
$T$. 
Invoking now Lemma~\ref{l.20} we conclude that 
the $j$-dimensional parallelepiped~$P(H)$ generated by $j$ directing segments 
of the face~$H$ is contained in~$T$. By Corollary~\ref{c.30}, 
for every directing segment~$[\bv, \ba]$ of $L$, 
the set~$T+ (\ba - \bv)$ is a tile from~$\cT$ and  hence, 
this set  lies in~$K$. Therefore, $P(H)+ (\ba - \bv) \, \subset \, K$. 
Thus, for every facet~$H$ of $L$ and  for every 
directing point $\ba$ of~$L$, the shifted parallelepiped 
$P(H)+ (\ba - \bv)$ lies in~$K$. 
Now we consider two cases. 

\smallskip 

1) $L$ is a simple cone. Then the $(j+1)$ shifted parallelepipeds~$P(H) \, + \, (\ba - \bv)$, 
where $H$ is a facet of~$L$ and $\ba\notin H$ is a directing point,  
form the boundary of the $(j+1)$-dimensional parallelepiped spanned by the 
directing segments of~$L$. This boundary lies in~$K$
and must intersect 
the open cone $S$. Hence $K\cap S \ne \emptyset$, which contradicts the definition of~$S$.
Thus, $S = \emptyset$ and so $L \subset K$.  
\smallskip 

2) $L$ is not  simple, i.e., it has at least $j+2$ edges. Let us show that 
this case is impossible. Applying Lemma~\ref{l.30} to the 
cone~$L$, we find two its $j$-dimensional faces $A, B$ and  two its 
different edges~$a, b$ with the directing points $\ba\in a, \bb \in b$ 
such  that the parallelepipeds 
$P(A) \, + \, (\ba - \bv)$ and $P(B) \, + \, (\bb - \bv)$
have a common interior point. The first one is contained 
in the tile~$T + (\ba - \bv)$, the second one is in the tile~$T + (\bb - \bv)$. 
Moreover, since $P(A)$ has a  layer in~$T$ (Lemma~\ref{l.10}), it follows that 
those two tiles have a common 
interior point. Hence, they coincide, which is impossible,  
because $(\ba - \bv) \ne (\bb - \bv)$, since those vectors are in 
different edges of~$L$.  

This completes the proof of the inductive step, which proves the proposition. 

{\hfill $\Box$}
\medskip 

Thus, we have proved that for every stationary vertex of~$G$, the corresponding 
corner of~$G$ is convex and simple. Since $G$ has at least one convex corner, it follows that it possesses  a simple convex corner. Now we can apply 
the main result of~\cite[Theorem 1.9]{YZ1}:  if a set  admits a self-affine tiling and has at least one convex polyhedral corner, then it is affinely similar to a union of integer shifts of the 
unit cube.  The next step is to show that a self-affine tiling of this set 
contains a tiling of a parallelepiped. 
\begin{lemma}\label{l.40}
Let a set $G$ be a union of integer shifts of a unit cube. 
Suppose~$G$ has a convex simple corner~$K$
at its extreme vertex and has a self-affine tiling~$\cT$;  
then there is an iteration of the tiling~$\cT$ whose subset 
 forms a tiling of some parallelepiped with the same corner~$K$.  
\end{lemma}
{\tt Proof} is in the Appendix.   
\begin{remark}\label{r.100}
{\em If among the unit cubes composing the set $G$ 
there is at least one separated from others, then 
Lemma~\ref{l.40} is obvious. Indeed, if an iteration of~$\cT$
has tiles of diameter smaller than one, then all tiles intersecting the separated cube 
do not intersect the others. Hence, they form a tiling of that cube.  
However, if there are no separated cubes in~$G$, then Lemma~\ref{l.40} is less obvious. 
}
\end{remark}

\medskip 

{\tt Proof of 
 Theorem~\ref{th.10}.} Let $G$ be a polyhedral set and $\cT$ be its self-affine tiling. 
 By Proposition~\ref{p.45}, there is some power of this tiling for which  
 $G$ has a stationary vertex~$\bv$. Proposition~\ref{p.50} asserts that the corner 
 $K$ at~$\bv$ is convex and simple. Hence, by \cite[Theorem~1.9]{YZ1}, 
 the set $G$ is equivalent to a union of several integer shifts of the unit cube. 
 All corners of this set are rectangular, i.e., equal to the cone~$\re^d_+$. 
 Now we invoke Lemma~\ref{l.40} and conclude that there is a power of the tiling~$\cT$
 whose subset forms a tiling of a parallelepiped~$P$ with the corner~$K$. 
 Hence, $P$ is a rectangular parallelepiped. Since the tile~$T$ covering the vertex~$\bv$
is affinely similar to $G$ and  
 has the same corner~$K$, it follows that after an affine transform with a diagonal matrix, 
 $T$ becomes a union of integer shifts of a unit cube, and several integer shifts 
 of $T$ cover a parallelepiped. Replacing each unit cube by its center we obtain a
 discrete tiling of a parallelepiped in~$\z^d$. According to the main result of \cite{Nath}, every discrete tiling is a direct product 
 of univariate discrete tilings of a segment of integer numbers. 
  Therefore, $\cT$
 is a direct product of $d$ tilings of a segment. 
Applying the classification of univariate tilings of a segment \cite[Theorem 8]{Zai} completes the proof.  

{\hfill $\Box$}
\medskip

\bigskip

\begin{center}
\large{\textbf{6. Proof of Theorem~\ref{th.20}}}
\end{center}
\medskip 

Let a polyhedral set~$G$ admit a self-affine tiling~$\cT$.  
 If $G$ is an integer attractor, then~$\cT \, = \, \{M^{-1}(G+\bss), \ \bss \in D\}$,  
$M$ is an integer expanding matrix, $D$ 
is  a set of digits, i.e., of representatives of  
quotient classes~$\z^d/M\z^d$.
\smallskip 

{\tt Proof of Theorem~\ref{th.20}}. 
Applying Theorem~\ref{th.10} we obtain that $G$
 is a direct product of~$d$ sets of the form~(\ref{eq.one-dim}). Hence, 
 $G$ is equivalent to a union of disjoint integer shifts of the unit cube. 
 For every~$n \ge 2$, the $n$th iteration of the tiling~$\cT$
 is defined by the matrix~$M^n$ and by the set of digits 
 $D_n = D + MD + \cdots + M^{n-1}D$. Clearly, all the elements of~$D_n$
 are from different quotient classes of~$\z^d/M^n\z^d$. 
 If $n$ is large enough, then the diameter of the set $M^{-n}G$ is smaller that one. 
 Hence, each  tile is contained in one of the unit cubes composing~$G$. 
 
 If $G$ contains more than one unit cube, then we take two of them~$C$
 and $C'$. We have 
  $C' = C+\ba$, where $\ba \in \z^d$.  Let $\bu$ be an arbitrary vertex of $C$ and 
 $T\in \cT^n$ be a tile containing~$\bu$. In the partial order defined by the corner 
 of the cube~$C$~at~$\bu$ the point~$\bu$ is the minimal point of~$C$. 
 Since $T \subset C$ it follows that $\bu$ is the minimal point of~$T$. 
 The vertex  $\bu+\ba$ of the cube~$C'$ is covered by another tile~$T'\in \cT^n$. 
 Let $T' = T + \bb, \ \bb \in \z^d$. By the same argument we show that 
 $\bu + \ba$ is the minimal point of~$T'$ with respect to the same order. 
 Since the parallel translation respects minimal points, 
 we have $\bu + \ba \, = \, \bu + \bb$  and so $\ba = \bb$. 
 Thus, $T' = T+\ba$. 
 Since~$\cT^n\,= \, \{M^{-n}(G + \bss), \ \bss \in D_n\}$, we  
 have $T\, = \, M^{-n}(G + \bss_1), \, T'\, = \, M^{-n}(G + \bss_1)$, 
 where $\bss_1, \bss_2 \in D_n$. 
 Thus, $T' \, = \, T\, +\, \ba\, =  \, T\, + \, M^{-n} (\bss_2 - \bss_1)$. 
 This means that $\ba\, =  \, M^{-n} (\bss_2 - \bss_1)$
and so $\bss_2 - \bss_1 \, = \, M^n\ba$, which is impossible since 
$\bss_1$ and $\bss_2$ are from different quotient classes of~$\z^{d}/ M^n\z^d$.

{\hfill $\Box$}
\medskip

 \bigskip 
 
\begin{center}
\large{\textbf{Acknowledgements}}
\end{center}
\bigskip 
The  work of  the  second  author  is  supported  by  the  Foundation  for  Advancement  of Theoretical Physics and Mathematics ``BASIS''. 
 
 \bigskip 
 
\begin{center}
\large{\textbf{Appendix}}
\end{center}
\bigskip 

{\tt Proof of Lemma~\ref{l.10}.} Let $\bz$ be an arbitrary point of~$G_L$ and an element $T' = T+\ba$
contain this point. Let a hyperplane $\{\bx \ | \ (\bc, \bx) = 0\}$
strongly separate~$L$ from $K$. Then it strongly separates $G_L$ from $K$. 
This means that  $(\bc, \bz) = (\bc, \bv) = 0$, while 
$(\bc, \bx) < 0$ for all points $\bx \in G\setminus G_L$. 
We have $\bz - \ba \in T$, hence $(\bc, \bz) - (\bc, \ba) \le 0$. 
Thus, $(\bc, \ba) \ge (\bc, \bz) = 0$, which implies  $(\bc, \bv + \ba) = 
(\bc , \ba) \ge 0$. 
However, $\bv + \ba \in T'$ and therefore, $\bv + \ba \in G$. 
This yields $(\bc, \bv + \ba) \le 0$ and consequently, $(\bc, \bv + \ba) = 0$. 
Thus, $\bv + \ba \in G_L$ and hence, 
the translation vector $\ba$ is parallel to~$L$. This is true for all 
tiles intersecting $G_L$. 

Thus, the sets $T'\cap L$, where 
$T'$ is a translation of~$T$ parallel to~$L$, are translations of~$T\cap L$. 
Clearly, they cover $G_L$.  
It remains to show that the interiors  (in~$L$) of those sets  are disjoint, 
from which it follows that they form a tiling of~$G_L$. To this end we prove 
that each interior point~$\bx \in T\cap L$ is contained in a neighbourhood 
$\cU(\bx) \subset G$. This will imply that the interiors of two different 
sets $T+\ba$ and $T+\bb$, where $\ba, \bb$ are parallel to $L$, do not intersect. 
Assume the contrary: the set of interior points from~$\bx \in T\cap L$ that 
do not have such a neighbourhood is nonempty. Denote this set by~$X$. 
Observe that $X$ does not intersect a neighbourhood of the vertex~$\bv$. 
Consider an arbitrary vector $\bh \in {\rm int}\, (L^*)$, where 
$L^* = \{\by \in {\rm span}\, (L) \ | \ (\by, \bx) \ge 0, \, \bx \in L\}$ 
is the dual cone in the linear span of~$L$. Since $X$ does not contain a neighbourhood 
of~$\bv$, we have $p = \inf_{\bx \in X} (\bh, \bx) > 0$. 
Take a positive number~$\varepsilon > 0$ which is less than $p$ and less than 
all numbers $(\bh, \ba), \ \ba \in L,\ T+\ba \in \cT$. Since the set
$\cT$ is finite, such a number exists. Now choose arbitrary~$\bx \in X$
such that $(\bh , \bx) < p + \varepsilon$. Since $\bx$ does not have a 
neighbourhood in~$G$, it must be covered by another  tile $T + \ba$, for which 
$\bx$ also does not have a neighbourhood in~$G$. Hence, 
$\bx - \ba \in X$, which is  impossible, since  $(\bh, \bx - \ba) < 
p + \varepsilon - \varepsilon = p$.   
 
 Since the affine similarity transform~$A$ respects~$L$ it follows that 
 all the sets~$T_i\cap G_L$, which are nonempty, are similar to~$G_L$ by the transform~$A$. 
 So, this tiling is self-affine. Finally, for every 
 $T'\in \cT, \, T'\ne T$, there is a bounded cone 
 in~$G$ with the base~$T\cap L$ that does not intersect~$T'$. Since the set~$\cT$ is finite, the intersection of those cones contains a bounded cone with the base~$T\cap L$. 
 This cone does not intersect other tiles from~$\cT$, hence, it lies in~$T$. 
 Therefore, this is a layer of the facet~$T\cap L$ in~$T$.

{\hfill $\Box$}
\medskip

{\tt Proof of Lemma~\ref{l.30}}. If $C$ is not simple, then there are 
two its  facets $A, B$ without a common $(d-2)$-dimensional face. 
Then there are two edges of the cone~$C$:  $a \subset A$ and $b \subset B$ 
such that the $(d-1)$ dimensional cones $C_a = {\rm co}\{a, B\}$ and 
$C_b = {\rm co}\{b, A\}$ have a common interior point. Clearly, $a$ is not in~$B$, 
otherwise $C_a  = B$ and it cannot have  common interior points 
with  $C_b$ since $A$ and $B$ have no common facets. 
Similarly, $b$ is not in~$A$. 
Since $a$ and $b$ are both in $C_a\cap C_b$ it follows 
that  the set $C_a\cap C_b$ has interior points arbitrary close to 
${\rm co}\{a, b\}$. Now denote by $\ba_i\, , \ i = 1, \ldots d-1$, 
the directing points on  the edges of the face~$A$, $\ba_1 \in a$. 
Analogously,   $\bb_i\, , \ i = 1, \ldots d-1$, 
are the directing points on  the edges of~$B$, $\bb_1 \in b$.
A common interior point $\bx$ of $C_a$ and $C_b$ is expressed as follows: 
\begin{equation}\label{eq.ab1}
\bx \quad = \quad \alpha\, \ba_1 \, + \, \sum_{i=1}^{d-1}\, t_i \bb_i \quad = \quad 
\beta\, \bb_1 \, + \, \sum_{i=1}^{d-1}\, s_i \ba_i\, , 
\end{equation}
 where all the  coefficients~$\alpha, \beta, t_i, s_i$ are strictly positive. 
 After multiplication of this equality by a positive constant it 
 can be assumed that $\alpha < 1$ and $\beta< 1$. 
 Moreover, choosing $\bx$ sufficiently close to~${\rm co}\{a, b\}$ we 
 may assume that all $t_i$ and $s_i$ are small, in particular, all of them
 are less than~$1$ and  
 $t_1 < \beta,\  s_1 < \alpha$.  Rewriting~(\ref{eq.ab1}) we obtain 
 $$
-  (\beta - t_1)\, \bb_1 \, + \, \sum_{i=2}^{d-1}\, t_i \bb_i \quad = \quad 
- (\alpha - s_1)\, \ba_1 \, + \, \sum_{i=2}^{d-1}\, s_i \ba_i\, . 
 $$
 Now we add $\ba_1 + \bb_1$ to both sides of this equality and get 
 \begin{equation}\label{eq.ab2}
\ba_1\, + \, \Bigl( 1 -  (\beta - t_1)\Bigr)\, \bb_1 \, + \, \sum_{i=2}^{d-1}\, t_i \bb_i \quad = \quad 
\bb_1 \, + \, \Bigl( 1 - (\alpha - s_1)\Bigr)\, \ba_1 \, + \, \sum_{i=2}^{d-1}\, s_i \ba_i\, . 
\end{equation}
 Note that the point $\Bigl( 1 -  (\beta - t_1)\Bigr)\, \bb_1 \, + \, \sum_{i=2}^{d-1}\, t_i \bb_i$ belongs to the interior of the parallelepiped~$P(B)$. 
 Indeed, this is a linear combination  of the vectors  $\bb_1, \ldots , \bb_{d-1}$
 and all coefficients of this combination are from the interval~$(0,1)$. 
 Hence, the left-hand side of equality~(\ref{eq.ab2}) is an interior 
 point of the parallelepiped~$\ba_1 + P(B)$. Analogously, the left-hand side  is an interior 
 point of~$\bb_1 + P(A)$. Thus, the parallelepipeds
 $\ba_1 + P(B)$ and $\bb_1 + P(A)$ possess a common interior point. 
 It remains to denote~$\ba = \ba_1, \, \bb= \bb_1$, which completes the proof.

{\hfill $\Box$}
\medskip

{\tt Proof of Lemma~\ref{l.40}}. By iterating the tiling, we may assume that 
$M$ is a diagonal matrix and the diameter of~$T$ is less than one. Let~$P_0$ be the 
parallelepiped spanned by the edges $[\bv, \bb_s] \, , \ s = 1, \ldots , d$
from Lemma~\ref{l.20}. Let us show that the parallelepiped~$P = M P_0$
is one we are looking for.   Since $K$ is simple, we identify it with~$\re^d_+$. 
The image of an integer shift of the unit cube under the action 
of~$M^{-1}$ will be called {\em brick}. Since $G$ consists of 
shifts of the unit cube, it follows that $T$ consists of bricks. 

We have~$T+\bb_i \in \cT$, therefore, the interior of the parallelepiped 
$P_0+ \bb_i $, being a part of the tile~$T+\bb_i$, does not intersect~$T$. 
Denote by $\{\be_j\}_{j=1}^d$  the canonical basis of~$\re^d$ and by 
$h_j$ the lengths of the edges of~$P$. Applying the similarity of $G$ and $T$
we conclude that the interior of the parallelepiped $P\, +\, h_1\be_i$ 
does not intersect~$G$. Among all unit cubes forming~$G$
we choose the ``highest'' ones with respect to the $i$th coordinate 
(whose center has the largest $i$th coordinate among all cubes). 
There are two possible cases: 

Suppose  none of  the highest cubes intersect the axis~$Ox_i$;  
then the tile $T'\in \cT$ containing the vertex~$h_i\be_i$ of~$P$
possesses a brick higher than that vertex, i.e., the $i$th coordinate of the center of that 
brick exceeds~$h_i$. Since the diameter of~$T'$ is smaller than one, that 
brick is contained in~$P+h_i\be_i$, which is impossible. 
Therefore,  for every~$i$, there exists the highest with respect to the $i$th 
coordinate cube in $G$ intersecting the axis~$Ox_i$. Hence, $T$ 
also has the highest brick (denote it by~$B_i$) which intersects 
the axis~$Ox_i$. 

If some tile~$T+\ba \in \cT$
intersecting the interior of~$P$ has a point whose $i$th coordinate 
exceeds~$h_i$, then the brick~$B_i + \ba$ is above the level $x_i = h_i$. 
On the other hand, all other coordinates of that brick are on the segments $[0, h_k], \, 
k\ne i$. Hence, $B_i + \ba \subset P+h_i \be_i$. Therefore, 
the brick $B_i + \ba $ is out of~$G$, which is impossible, since this is a part of the 
tile~$T+\ba$. Thus, all points of each tile intersecting~$P$
have the $i$th coordinate at most~$h_i$. Applying this argument for  
all~$i$, we see that every tile intersecting $P$ lies in~$P$. 
Hence, those tiles form a tiling of~$P$.

{\hfill $\Box$}

\end{document}